\documentclass[10pt]{amsart}
\usepackage[all]{xy}

\usepackage{amssymb, amsthm,amsmath,dsfont}
\usepackage{enumerate}
\usepackage{stmaryrd}
\usepackage{mathrsfs}

\begin{document}

\theoremstyle{plain}
\newtheorem{theorem}{Theorem}
\newtheorem{coro}[theorem]{Corollary}
\newtheorem{lemma}[theorem]{Lemma}
\newtheorem{fact}[theorem]{Fact}
\newtheorem{proposition}[theorem]{Proposition}
\newtheorem{ques}[theorem]{Question}
\newtheorem{step}{Claim}
\newtheorem{observation}[theorem]{Observation}
\newtheorem{conj}[theorem]{Conjecture}

\theoremstyle{remark}
\newtheorem{example}[theorem]{Example}
\newtheorem{remark}[theorem]{Remark}
\newtheorem*{notation}{Notation}

\numberwithin{theorem}{section}
\numberwithin{equation}{section}

\theoremstyle{definition}
\newtheorem{definition}[theorem]{Definition}

\newcommand{\F}{\mathbb{F}}
\newcommand{\Z}{\mathbb{Z}}
\newcommand{\N}{\mathbb{N}}
\newcommand{\Q}{\mathbb{Q}}
\newcommand{\K}{\mathbb{K}}
\newcommand{\sB}{\mathcal{B}}
\newcommand{\sS}{\mathcal{S}}

\newcommand{\sep}{\mathrm{sep}}
\newcommand{\chr}{\mathrm{char}}
\newcommand{\ext}{\mathrm{Ext}}
\newcommand{\Span}{\mathrm{Span}}
\newcommand{\Hom}{\mathrm{Hom}}
\newcommand{\Gal}{\mathrm{Gal}}

\newcommand{\dbl}{[\![}
\newcommand{\dbr}{]\!]}
\newcommand{\dbll}{(\!(}
\newcommand{\dblr}{)\!)}
\newcommand{\dbml}{\langle\!\langle}
\newcommand{\dbmr}{\rangle\!\rangle}
\newcommand{\gr}{{\rm gr}}
\newcommand{\res}{{\rm res}}
\newcommand{\sU}{{\mathcal U}}
\newcommand{\FpG}{\F_p[G]}
\newcommand{\FppG}{\F_p\dbl G\dbr}
\newcommand{\FpX}{\F_p\langle X\rangle}
\newcommand{\FppX}{\F_p\dbml X\dbmr}

\newcommand{\kernel}{\mathrm{ker}}
\newcommand{\image}{\mathrm{im}}
\newcommand{\iid}{\mathrm{id}}
\newcommand{\rk}{\mathrm{rk}}
\newcommand{\dd}{\mathrm{d}}


\newcommand{\hac}{\hat c}
\newcommand{\hatheta}{\hat\theta}

\title[One-relator pro-$p$ Galois groups and the koszulity conjectures]
 {One-relator maximal pro-$p$ Galois groups \\ and the koszulity conjectures}
 
\author{Claudio Quadrelli}

\dedicatory{To John P. Labute, with respect and admiration.}

\address{Department of Mathematics and Applications, University of Milano Bicocca, Via R. Cozzi 55 - Ed. U5, 20125 Milan, Italy EU}
\email{claudio.quadrelli@unimib.it}
\date{\today}

\begin{abstract}
Let $p$ be a prime number and let ${\K}$ be a field containing a root of 1 of order $p$.
If the absolute Galois group $G_{\K}$ satisfies $\dim H^1(G_{\K},\F_p)<\infty$ and $\dim H^2(G_{\K},\F_p)=1$,
we show that L.~Positselski's and T.~Weigel's Koszulity conjectures are true for ${\K}$.
Also, under the above hypothesis we show that the $\F_p$-cohomology algebra of
$G_{\K}$ is the quadratic dual of the graded algebra $\gr_\bullet\F_p[G_{\K}]$,
induced by the powers of the augmentation ideal of the group algebra $\F_p[G_{\K}]$,
and these two algebras decompose as products of elementary quadratic algebras.
Finally, we propose a refinement of the Koszulity conjectures, analogous to 
I.~Efrat's Elementary Type Conjecture.
\end{abstract}




\maketitle

\section{Introduction}
\label{sec:intro}

Let $A_\bullet:=\bigoplus_{n\geq0}A_n$ be a non-negatively graded algebra of finite type over a field $\Bbbk$.
The algebra $A_\bullet$ is called {\sl quadratic} if it is generated in degree 1 (i.e., every element is a combination
of products of elements of $A_1$) and its defining relations are homogeneous relations of degree 2.
Namely, one may write $A_\bullet\simeq T_\bullet(A_1)/(\Omega)$, where $T_\bullet(A_1)$ denotes the tensor
algebra generated by $A_1$ and $(\Omega)$ is the two-sided ideal generated by a subset $\Omega\subseteq A_1^{\otimes2}$.
A quadratic algebra $A_\bullet$ comes equipped with its {\sl quadratic dual} $A_\bullet^!$, which is the quadratic
algebra generated by the dual $A_1^*$ of $A_1$, and with defining relations the orthogonal complement
$\Omega^\perp\leq(A_1^*)^{\otimes2}$ of $\Omega$ (cf. \cite[\S~1.2]{PoliPosi}).
 
Quadratic algebras gained great importance in Galois theory after the proof of the celebrated Bloch-Kato conjecture by M.~Rost 
and V.~Voevodsky (see \cite{rost,BK,weibel,weibel:art}): from it one deduces that, given a prime number $p$ and a field $\K$
containing a root of 1 of order $p$, the $\F_p$-cohomology algebra $H^\bullet(G_{\K})=\bigoplus_{n\geq0}H^n(G_{\K},\F_p)$
of the absolute Galois group $G_{\K}$ of $\K$, endowed with the (graded-commutative) cup-product
\[
 \cup \colon H^r(G_{\K},\F_p)\otimes H^s(G_{\K},\F_p)\longrightarrow H^{r+s}(G_{\K},\F_p),
\]
is a quadratic algebra over the finite field $\F_p$.
This led to the achievement of new results on the structure of maximal pro-$p$ Galois groups of fields
(see, e.g., \cite{KLM,cem,cq:bk,CMQ,cq:qc}).
Understanging which pro-$p$ groups may occur as maximal pro-$p$ Galois groups of fields
is one of the main open problems in modern Galois theory.

The class of {\sl Koszul algebras} is a very peculiar class of quadratic algebras, singled out by S.B.~Priddy in \cite{priddy}.
A quadratic algebra $A_\bullet$ is called Koszul if $\Bbbk$, as trivial graded $A_\bullet$ module concentrated in degree 0, has a linear resolution (see Definition~\ref{defi:koszul} for the formal definition).
Koszul algebras arise in various fields of mathematics, and they have an uncommonly nice cohomological
behavior --- e.g., the cohomology of a Koszul algebra is just its quadratic dual.
For the formal definition and properties of Koszul algebras we direct the reader to \cite[Ch.~2]{PoliPosi}
and \cite[\S~2]{MQRTW}.
Koszul algebras appeared on the Galois stage thanks to the work of L.~Positselski and A.~Vishik
(see, e.g., \cite{posivis:koszul,positselsky:koszul1}).
In particular, in \cite{positselsky:koszul2} Positselski conjectured that $H^\bullet(G_{\K})$ is Koszul,
for $\K$ containing a root of 1 of order $p$, and he proved that this is the case if $\K$ is a number field.
In view of Positselski's conjecture and of the cohomological properties of Koszul algebras, it is natural to ask
what the quadratic dual of $H^\bullet(G_{\K})$ may look like.

For a profinite group $G$, let $$\gr_\bullet\FpG:=\bigoplus_{n\geq0}I^n/I^{n+1}$$ denote the
{\sl graded group algebra} of $G$, with $I\trianglelefteq\FpG$ the augmentation ideal and $I^0=\FpG$.
If $\K$ contains a root of unity of order $p$ and the maximal pro-$p$ Galois group of $\K$ is finitely
generated, $H^\bullet(G_{\K})$ and $\gr_\bullet\F_p[G_{\K}]$ are expected to be quadratically dual to each other,
and Th.~Weigel conjectured in \cite{thomas:proc} that also $\gr_\bullet\F_p[G_{\K}]$ is Koszul.
Altogether, one has the following three ``Koszulity conjectures''.

\begin{conj}\label{conj}
 Let $\K$ be a field containing a root of 1 of order $p$ with $H^1(G_{\K},\F_p)$ finite.
\begin{itemize}
  \item[(i)]  The $\F_p$-cohomology algebra $H^\bullet(G_{\K})$ is Koszul (cf. \cite{positselsky:koszul2}). 
 \item[(ii)]  The graded group algebra $\gr_\bullet\F_p[G_{\K}]$ is Koszul (cf. \cite{thomas:proc}).
\item[(iii)]  $H^\bullet(G_{\K})^!\simeq\gr_\bullet\F_p[G_{\K}]$ (cf. \cite{MQRTW}).
\end{itemize}
\end{conj}

In \cite{MQRTW}, it is shown that the above conjecture holds true for the class of pro-$p$ groups {\sl of elementary type}
(introduced by I.~Efrat in \cite[\S~3]{ido:ETC}), which includes Demushkin pro-$p$ groups,
and which is particularly significan from a Galois-theoretic point of view.

Following the trail drawn in \cite{MQRTW}, in this paper we study the $\F_p$-cohomology algebra and
the graded group algebra for the class of finitely generated {\sl one-relator} pro-$p$ groups with quadratic $\F_p$-cohomology.
A pro-$p$ group $G$ is said to be one-relator if it has a minimal presentation with a single defining relation.
The fundamental example we keep in mind while dealing with one-relator pro-$p$ groups is the class of Demushkin
groups --- introduced by S.P.~Demushkin and calssified completely by J.-P.~Serre and J.P.~Labute ---,
which arise as maximal pro-$p$ Galois groups of local fields (see, e.g., \cite{labute:demushkin}
and \cite[\S~III.9]{nsw:cohm}).
Our investigation on one-relator pro-$p$ groups aims at proving the following.

\begin{theorem}\label{thm:main}
Let $\K$ be a field containing a root of 1 of order $p$ such that $H^1(G_{\K},\F_p)$ is finite and
$\dim H^2(G_{\K},\F_p)=1$.
Then there exists a closed subgroup $S$ of $G_{\K}$ with $\mathrm{cd}_p(S)=1$,
generating a closed normal subgroup $N\trianglelefteq G_{\K}$, such that:
\begin{itemize}
 \item[(i)] the maximal pro-$p$ quotient of $G_{\K}/N$ is a Demushkin group;
 \item[(ii)] $H^\bullet(G_{\K})\simeq H^\bullet(S)\sqcap H^\bullet(G_{\K}/N)$;
 \item[(iii)] $\gr_\bullet\F_p[G_{\K}]\simeq\gr_\bullet\F_p[S]\sqcup \gr_\bullet\F_p[G_{\K}/N]$.
Moreover, also $\mathrm{cd}_p(N)=1$, unless $p=2$ and $\sqrt{-1}\notin \K$. 
\end{itemize}
\end{theorem}

Here $\sqcap$ and $\sqcup$ denote the direct product and free product in the category of quadratic algebras (see \S~2.2).
The assumption on $H^n(G_{\K})$, $n=1,2$, amounts to say that the maximal pro-$p$ quotient of $G_{\K}$
is a finitely generated one-relator pro-$p$ group.
Theorem~\ref{thm:main} may be considered as an ``algebras analogue'' of a result by T. W\"urfel
on one-relator pro-$p$ groups which occur as absolute Galois groups (cf. \cite{wurf}).

More importantly, Theorem~\ref{thm:main} implies that the ``Koszulity conjectures''
find positive solution under the assumption that $H^1(G_{\K},\F_p)$ is finite and
$H^2(G_{\K},\F_p)$ has dimension 1. 

\begin{coro}\label{cor:main}
Let $\K$ be as in Theorem~\ref{thm:main}. Then Conjecture~\ref{conj} holds true for $\K$. 
\end{coro}

Finally, in analogy with I. Efrat’s Elementary Type Conjecture (cf. \cite{ido:ETC}), we define the classes
of Koszul algebras of {\sl $H$-elementary type} and of {\sl $G$-elementary type} as the minimal classes of Koszul algebras
generated by some ``basic'' quadratic algebras via elementary products (cf. Definition~\ref{defi:pbwetc}).
In view of the results obtained in \cite{MQRTW} and in the present paper, we ask the following.

\begin{ques}
Let $\K$ be a field containing a root of unity of order $p$ such that $H^1(G_{\K} ,\F_p)$ is finite.
\begin{itemize}
 \item[(i)] Is the $\F_p$-cohomology $H^\bullet(G_{\K})$ a Koszul algebra of $H$-elementary type?
 \item[(ii)] Is the graded group algebra $\gr_\bullet\F_p[G_{\K}]$ a Koszul algebra of $G$-elementary type?
\end{itemize}
\end{ques}

\subsection*{Acknowledgment}
{\small The author wishes to express his gratitude toward the anonymous referee for her/his careful work and useful comments, and toward
J.~Min\'a\v c, F.W. Pasini, N.D. T\^an and Th. Weigel, as his interest for Koszulity in Galois
cohomology originates from the collaboration with them.
Moreover, the author thanks I. Efrat and J.P. Labute for their precious suggestions, and P. Guillot, C. Maire, E. Matzri, D. Neftin, M. Schein and I. Snopce for their interest.}


\section{Quadratic algebras}

Throughout the paper every graded algebra $A_\bullet=\bigoplus_{n\in\Z}A_n$ is assumed to be
unitary associative over the finite field $\F_p$, and non-negatively graded of finite-type,
i.e., $A_0=\F_p$, $A_n=0$ for $n<0$ and $\dim(A_n)<\infty$ for $n\geq1$.

In this section we give only the most basic definitions and results on quadratic algebras and Koszul algebras,
and some examples, which will be sufficient for our investigation.
For a complete account on cohomology of graded algebras and Koszul algebras, we refer to \cite[\S~2]{MQRTW},
and also to the first chapters of the books \cite{PoliPosi} and \cite{lodval}.


\subsection{Quadratic algebras and quadratic duals}

Given a vector space $V$ of finite dimension, let $T_\bullet(V)$ denote
the graded tensor algebra generated by $V$, endowed with the multiplication induced by the tensor product.
Moreover, let $V^*=\Hom_{\F_p}(V,\F_p)$ denote the dual space of $V$.
Since $\dim V<\infty$, one may identify $(V\otimes V)^*=V^*\otimes V^*$.

\begin{definition}\label{defi:algebras} 
A graded algebra $A_\bullet=\bigoplus_nA_n$ is said to be {\sl quadratic} if $A_\bullet$ is isomorphic to the quotient 
$T_\bullet(V)/(\Omega)$ for some finitely generated vector space $V$ isomorphic to $A_1$, and some subset $\Omega\subseteq V\otimes V$,
with $(\Omega)\trianglelefteq T_\bullet(V)$ the two-sided ideal generated by $\Omega$.
We write $A_\bullet=Q(V,\Omega)$.

For a quadratic algebra $A_\bullet=Q(V,\Omega)$, let $\Omega^\perp\subseteq (V\otimes V)^*$ denote the orthogonal space
of $\Omega$, i.e., $\Omega^\perp=\{\alpha\in(V\otimes V)^*\mid \alpha(w)=0\text{ for all }w\in \Omega\}$.
The {\sl quadratic dual} of $A_\bullet$, denoted by $A_\bullet^!$, is the quadratic algebra
$Q(V^*,\Omega^\perp)$.
\end{definition}

Note that for a quadratic algebra $A_\bullet$ one has $(A_\bullet^!)^!=A_\bullet$.

\begin{definition}\label{defi:koszul}
 For a graded algebra $A_\bullet=\bigoplus_{n\geq0}A_n$, the $\F_p$-cohomology is defined as the direct sum of the derived functors if the functor $\mathrm{Hom}_{A_\bullet}(\_,\F_p)$ evaluated on $\F_p$.
 Namely, it is the {\sl bigraded} algebra
 \[
  \bigoplus_{i,j}\mathrm{Ext}_{A_\bullet}^{ij}(\F_p,\F_p),\qquad i,j\geq0,
 \]
where one grading is induced by the grading of $A_\bullet$, and the other grading is the cohomological grading (cf. \cite[\S~2.2]{MQRTW}).
The algebra $A_\bullet$ is {\sl Koszul} if the cohomology is concentrated on the diagonal, i.e., if $\mathrm{Ext}_{A_\bullet}^{ij}(\F_p,\F_p)=0$ for $i\neq j$.
\end{definition}

Every Koszul algebra is quadratic.
Moreover, a quadratic algebra $A_\bullet$ is Koszul if, and only if, the dual algebra $A_\bullet^!$ is Koszul 
(cf. \cite[\S~2.3]{MQRTW}).


\subsection{Examples and constructions}

Henceforth $V$ denotes a vector space of finite dimension $d$.

\begin{example}\label{ex:tensor} 
 The tensor algebra $T_\bullet(V)$ and the quadratic algebra $Q(V,V^{\otimes2})$, called the
 {\sl trivial} quadratic algebra, are Koszul, and $Q(V,V^{\otimes2})^!=T_\bullet(V^*)$,
 and conversely (cf. \cite[Examples~3.2.5]{lodval}).
\end{example}

Let $ X =\{X_1,\ldots,X_d\}$ be a set of indeterminates.
The free associative algebra $\FpX$ --- i.e., the algebra of polynomials on the
non commutative indeterminates $X$ --- comes endowed with the grading induced by the subspaces of homogeneous polynomials.
We may identify $ X $ with a fixed basis of $V$, and such identification induces an isomorphism
of quadratic algebras $T_\bullet(V)\simeq\FpX$.

\begin{example}\label{ex:symmetric exterior} 
 The symmetric algebra $S_\bullet(V)=Q(V,\Omega_S)$ and the exterior algebra $\Lambda_\bullet(V)=Q(V,\Omega_\Lambda)$,
 where \[\Omega_S=\{vw-wv\mid v,w\in V\}\qquad\text{and}\qquad\Omega_\Lambda=\{vw+wv\mid v,w\in V\},\]
 are Koszul, and $\Lambda_\bullet(V)^!=S_\bullet(V^*)$, and conversely (cf. \cite[Examples~3.4.12]{lodval}).
\end{example}

Given two quadratic algebras $A_\bullet=Q(A_1,\Omega_A)$ and $B_\bullet=Q(B_1,\Omega_B)$, one has 
the following contrsuctions (cf. \cite[Example~2.5]{MQRTW}).
 \begin{itemize}
 \item[(a)] The {\sl direct product} of $A_\bullet$ and $B_\bullet$ is the quadratic algebra
 $A_\bullet\sqcap B_\bullet=Q(A_1\oplus B_1,\Omega)$, 
 with $$\Omega=\Omega_A\cup\Omega_B\cup (A_1\otimes B_1)\cup (B_1\otimes A_1).$$
 \item[(b)] The {\sl free product} of $A_\bullet$ and $B_\bullet$ is the quadratic algebra
$A_\bullet\sqcup B_\bullet=Q(A_1\oplus B_1,\Omega_A\cup\Omega_B)$.
 \item[(c)] The {\sl symmetric tensor product} of $A_\bullet$ and $B_\bullet$ is the quadratic algebra
$A_\bullet\otimes B_\bullet=Q(A_1\oplus B_1,\Omega)$,
with $\Omega=\Omega_A\cup\Omega_B\cup\Omega_S$, where $$\Omega_S=\{ab-ba, a\in A_1,b\in B_1\}.$$
\item[(d)] The {\sl skew-symmetric tensor product} of $A_\bullet$ and $B_\bullet$ is the quadratic algebra
$A_\bullet\wedge B_\bullet=Q(A_1\oplus B_1,\Omega)$,
with $\Omega=\Omega_A\cup\Omega_B\cup\Omega_\wedge$, where $$\Omega_\wedge=\{ab+ba, a\in A_1,b\in B_1\}.$$
\end{itemize}

One has the following (cf. \cite[\S~3.1]{PoliPosi}).

\begin{proposition}\label{prop:operations koszul}
 If both $A_\bullet$ and $B_\bullet$ are Koszul, then also the algebras $A_\bullet\sqcap B_\bullet$,
 $A_\bullet\sqcup B_\bullet$, $A_\bullet\otimes B_\bullet$ and $A_\bullet\wedge B_\bullet$ are Koszul.
Moreover, one has
\begin{itemize}
 \item[(i)] $(A_\bullet\sqcap B_\bullet)^!\simeq A^!_\bullet\sqcup B_\bullet^!$ and 
$(A_\bullet\sqcup B_\bullet)^!\simeq A^!_\bullet\sqcap B_\bullet^!$,
 \item[(ii)] $(A_\bullet\otimes B_\bullet)^!\simeq A^!_\bullet\wedge B_\bullet^!$ and 
$(A_\bullet\wedge B_\bullet)^!\simeq A^!_\bullet\otimes B_\bullet^!$.
\end{itemize}
\end{proposition}


\section{Pro-$p$ groups and graded algebras}
\label{sec:3}

The aim of this section is to provide an effective way to describe the graded group algebra of the pro-$p$ groups
we are studying.

For a pro-$p$ group $G$ and $n\geq1$, $G^n$ denotes the closed subgroup of $G$ generated by $n$th powers,
and $[G,G]$ denotes the closed commutator subgroup of $G$.

\subsection{The $p$-Zassenhaus filtration}
The following is an extract from  \cite[\S~2]{jochen:massey}.

Let $F$ be a free finitely generated pro-$p$ group with basis $\mathcal{X}=\{x_1,\ldots,x_d\}$.
The complete group algebra $\F_p\dbl F\dbr$ is defined as $\F_p\dbl F\dbr=\varprojlim_U\F_p[F/U]$,
where $U$ runs through the open normal subgroups of $F$.
The assignment $x\mapsto x-1$ induces an embedding of sets $F\hookrightarrow \F_p\dbl F\dbr$ (this map is not a morphism of groups).

Let $\FppX$ denote the algebra of formal power series in the non-commuting indeterminates $X=\{X_1,\ldots,X_d\}$.
Then there is an isomorphism of topological algebras $\phi\colon \F_p\dbl F\dbr$, given by $x_i\mapsto 1+X_i$.
The composition of the embedding $F\hookrightarrow\F_p\dbl F\dbr$ with the isomorphism $\phi$ is the Magnus
embedding $\psi\colon F\hookrightarrow\FppX$.

Let $I(X)$ denote the augmentation ideal of $\FppX$, i.e., $I(X)$ is the two-sided ideal $(X_1,\ldots,X_d)$.
The $p$-Zassenhaus filtration of $F$ is the filtration given by the subgroups
\[
 F_{(n)}=\{x\in F\mid \phi(x-1)\in I(X)^n\}=\{x\in F\mid \psi(x)\in I(X)^n\},\qquad n\geq1.
\]
In particular, one has $F_{(1)}=G$, $F_{(2)}=F^p[F,F]$ (namely, $F_{(2)}$ is the Frattini subgroup of $F$), and 
\begin{equation}\label{eq:G3}
 F_{(3)}=\begin{cases} F^p[[F,F],F] & \text{if }p\neq2 \\ F^4[F,F]^2[[F,F],F] & \text{if }p=2.\end{cases}
\end{equation}
Moreover, the $p$-Zassenhaus filtration of $F$ is the fastest descending series of $F$ starting at $F$ and such that
\begin{equation}\label{eq:propzass}
 [F_{(n)},F_{(m)}]\subseteq F_{(n+m)}\qquad\text{and}\qquad F_{(n)}^p\subseteq F_{(np)}
\end{equation}
for every $n,m\geq1$ (cf. \cite[\S~11.1]{ddsms}), and the quotient $F_{(n)}/F_{(n+1)}$ is finite.

Moreover, one has a canonical isomorphism of graded algebras 
\[ \gr\FppX:=\bigoplus_{n\geq0}I(X)^n/I(X)^{n+1}\overset{\sim}{\longrightarrow} \FpX, \]
so that for a series $f\in\FppX$ such that $f\in I(X)^n\smallsetminus I(X)^{n+1}$, one may consider 
the class $f+I(X)^{n+1}$ as a homogeneous polynomial of $\FpX$ of degree $n$.

For an element $x\in F_{(n)}\smallsetminus  F_{(n+1)}$, the class $\psi(x)+I(X)^{n+1}$, considered as a
homogeneous polynomial of $\FpX$ of degree $n$, is called the {\sl initial form} of $x$ in $\FpX$.
Thus, one may consider the quotient $F_{(n)}/F_{(n+1)}$ as a subspace of the space of homogeneous polynomials
on $X$ of degree $n$.
In particular, the initial form of $x_i$ is $X_i$, and the initial form of the commutator $[x_i,x_j]=x_i^{-1}x_j^{-1}x_ix_j$
is the algebra commutator $-[X_i,X_j]=X_iX_j+X_jX_i$, for every $i,j\in\{1,\ldots,d\}$.


\subsection{Restricted lie algebras}
In this subsection we give a description of the graded group algebra of a finitely generated pro-$p$
group as quotient of the algebra $\FpX$.

For $X=\{X_1,\ldots,X_d\}$, one can make the algebra $\FpX$ into a Lie algebra, by setting Lie brackets
$[f_1,f_2]=f_1f_2-f_2f_1$ for $f_1,f_2\in\FpX$.

\begin{definition}\label{defi:restricted}
Let $L$ be a Lie subalgebra of $\FpX$.
Then $L$ is said to be {\sl restricted} if $f^p\in L$ for each element $f\in L$.
In particular, the restricted Lie algebra $L(X)$ is the restricted Lie subalgebra of $\FpX$ generated by $X$,
it is a {\sl free restricted Lie algebra}, and $\FpX$ is its {\sl universal envelope}.
\end{definition}

(For the general definition of restricted Lie algebra see \cite[\S~V.7]{jacobson} and \cite[\S~12.1]{ddsms}.)
Let $F$ be a free pro-$p$ group with basis $\mathcal{X}=\{x_1,\ldots,x_d\}$.
For every $n\geq1$, the subspace of $L(X)$ of the homogeneous elements of degree $n$ is the image of the quotient
$F_{(n)}/F_{(n+1)}$ via $\psi$, and we may identify $L(X)$ with $\bigoplus_{n\geq0}F_{(n)}/F_{(n+1)}$
(cf. \cite[Rem.~2.3]{jochen:massey}).

Let $G$ be a pro-$p$ group with presentation
 \begin{equation}\label{eq:presentation G}
 \xymatrix{ 1\ar[r] & R\ar[r]& F\ar[r] & G\ar[r] & 1 }.
\end{equation}
Recall that a presentation \eqref{eq:presentation G} is minimal when $R\subseteq F_{(2)}$ (cf., e.g., \cite[p.~215]{CMQ}).
A subset $\{r_1,\ldots,r_m\}$ of $F$ is said to be a set of defining relations for $G$ if it generates 
$R$ as closed normal subgroup of $F$.

Set 
\[
 I(R)=\bigoplus_{n\geq2}\dfrac{(R\cap F_{(n)})F_{(n+1)}}{F_{(n+1)}}.
\]
Since $R$ is a normal subgroup of $F$, $I(R)$ is a {\sl restricted ideal} of $L(X)$, i.e., it is an ideal in the sense of a Lie algebra, with the further condition that $f^p\in I(R)$ for each $f\in I(R)$ (cf. \cite[\S~7.2]{MQRTW}).
We set $L(G)=L(X)/I(R)$, and $G_{(n)}=F_{(n)}/(R\cap F_{(n)})$ for $n\geq1$, i.e., 
\[ L(G)=\bigoplus_{n\geq1} G_{(n)}/G_{(n+1)}=L(X)/I(R)\]
(for the original definition of the $p$-Zassenhaus filtration $\{G_{(n)}\}_{n\geq1}$ of a pro-$p$ group $G$ see \cite[\S~11.1]{ddsms}).
The graded group algebra $\gr_\bullet\F_p[G]$ and the restricted Lie algebra $L(G)$ related by Jennings' theorem
(cf. \cite[\S~12.2]{ddsms}, see also \cite[Thm.~3.9]{MQRTW})

\begin{proposition}\label{prop:Jennings}
 Let $G$ be as above.
 Then $\gr_\bullet\F_p[G]$ is isomorphic to the universal envelope of the restricted Lie algebra $L(G)$. 
\end{proposition}

Proposition~\ref{prop:Jennings} and \cite[Prop.~2.1]{jochen:massey} yield the following.

\begin{proposition}\label{prop:rest}
 Let $G$ be a finitely generated pro-$p$ group with presentation \eqref{eq:presentation G}.
Then one has a short exact sequence of graded $\F_p$-algebras 
\[
\xymatrix{ 0\ar[r] & \mathcal{I} \ar[r] & \FpX \ar[r] & \gr_\bullet\F_p[G]  \ar[r]& 0},
\]
where $\mathcal{I}\trianglelefteq\FpX$ is the two-sided ideal (as ideal of an associative algebra)
generated by $I(R)$ as subset of $\FpX$.
\end{proposition}


\subsection{Mild pro-$p$ groups}

Let $X=\{X_1,\ldots,X_d\}$ be a set of non-commuting indeterminates.
For $m\geq1$, let $\rho_1,\ldots,\rho_m\in\FpX$ be homogeneous polynomials of degree $\deg(\rho_i)=s_i\geq2$
for each $i$.
Let $\mathcal{I}\trianglelefteq\FpX$ denote the two-sided ideal generated by $\rho_1,\ldots,\rho_m$,
and set $A_\bullet=\FpX/\mathcal{I}$.
The sequence $\{\rho_1,\ldots,\rho_m\}$ is called {\sl strongly free} if one has an identity of formal power series
\[
 \sum_{n\geq0}\dim(A_n)\cdot t^n=\dfrac{1}{1-dt+(t^{s_1}+\ldots+t^{s_m})}
\]
(cf. \cite[Def.~2.6]{jochen:massey}).

\begin{definition} 
 Let $G$ be a finitely generated pro-$p$ group with minimal presentation \eqref{eq:presentation G}, where $F$
has a basis $\mathcal{X}=\{x_1,\ldots,x_d\}$.
The group $G$ is said to be {\sl mild} if there exists a finite set of defining relations $\{r_1,\ldots,r_m\}$
such that $\{\rho_1,\ldots,\rho_m\}$ is a strongly free sequence in $\FpX$, with $\rho_i$
the initial form of $r_i$ for each $i$.
Such a presentation is called mild.
\end{definition}

The following result is fundamental for dealing with the algebras arising from mild pro-$p$ groups ---
in particular, it provides an effective description for the ideal $\mathcal{I}$ as in Proposition~\ref{prop:rest}
(cf. \cite[Thm.~2.12]{jochen:massey}).

\begin{proposition}\label{prop:mild}
 Let $G$ be a mild pro-$p$ group, with mild presentation \eqref{eq:presentation G}.
\begin{itemize}
 \item[(i)] $G$ has cohomological dimension 2, i.e., $H^n(G,\F_p)=0$ for $n\geq3$.
 \item[(ii)] One has an isomorphism of graded algebras $\gr_\bullet\FpG\simeq\FpX/\mathcal{I}$,
where $X=\{X_1,\ldots,X_d\}$ and $\mathcal{I}$ is the two-sided ideal generated by the initial forms of 
the defining relations of $G$.
\end{itemize}
\end{proposition}

Note that for a generic pro-$p$ group $G$ the dieal $\mathcal{I}$ as in Proposition~\ref{prop:rest} may not be generated only by the initial forms of a set of defining relations.

\begin{example}\label{ex:RAAG} 
For $p$ odd let $G$ be the pro-$p$ group generated by $\{x_1,x_2,x_3,x_4\}$ and subject to the defining relations
 $$[x_1,x_2]x_1^{\alpha_{12}}x_2^{\beta_{12}}=[x_2,x_3]x_2^{\alpha_{23}}x_3^{\beta_{23}}
 =[x_3,x_4]x_3^{\alpha_{34}}x_4^{\beta_{34}}=[x_4,x_1]x_4^{\alpha_{41}}x_1^{\beta_{41}}=1,$$
for some $\alpha_{ij},\beta_{ij}\in p\Z_p$.
Tge group may be reslized as the {\sl generalised right-angled Artin pro-$p$ group} associated to the square graph with vertices $\{x_1,x_2,x_3,x_4\}$ (cf. \cite[Def.~5.3]{QSV}). 
By \cite[Thm.~F]{QSV}, $G$ is mild.
Since the initial form of the defining relations above are the commutators $[X_1,X_2],\ldots,[X_4,X_1]$, by Proposition~\ref{prop:mild} one has 
$$\gr_\bullet\FpG\simeq\frac{\FpX}{([X_1,X_2],[X_2,X_3],[X_3,X_4],[X_4,X_1])},$$
where $X=\{X_1,\ldots,X_4\}$.
Consequently, by \cite[\S~4.2.2]{thomas:FP} the algebra $\gr_\bullet\F_p[G]$ is Koszul.
\end{example}


\section{One-relator pro-$p$ groups}\label{sec:4}

For a pro-$p$ group $G$ we shall denote the $\F_p$-cohomology groups simply by $H^n(G)$ for every $n\geq0$.
Thus $H^0(G)=\F_p$, and given a presentation \eqref{eq:presentation G}, one has the following isomorphisms of vector spaces:
\begin{equation}\label{eq:H1H2}
 \begin{split}
  H^1(G) &\simeq H^1(F)\simeq(G/G_{(2)})^*, \\ H^2(G) &\simeq H^1(R)^G\simeq(R/R^p[R,F])^*    \end{split}  \end{equation}
(cf. \cite[Prop.~3.9.1 and Prop.~3.5.9]{nsw:cohm}).
The $\F_p$-cohomology of a pro-$p$ group comes endowed with the {\sl cup-product} 
\[
 \xymatrix{ H^i(G)\times H^j(G)\ar[r]^-{\cup} & H^{i+j}(G) }
\]
which is graded-commutative, i.e., $\alpha\cup\beta=(-1)^{ij}\beta\cup\alpha$ for $\alpha\in H^i(G)$ and $\beta\in H^j(G)$.
For further facts on cohomology of pro-$p$ groups we refer to \cite[\S~III.9]{nsw:cohm}.

Henceforth $G$ is assumed to be a finitely generated one-relator pro-$p$ group.
Thus, if \eqref{eq:presentation G} is a minimal presentation, then $R$ is generated by a single defining relation
$r\in F_{(2)}$.
Given a fixed basis $\mathcal{X}=\{x_1,\ldots,x_d\}$ of $F$, and a basis $\mathcal{B}=\{\chi_1,\ldots,\chi_d\}$ of $H^1(G)$ dual to $\mathcal{X}$ (i.e., $\chi_i(x_j)=\delta_{ij}$), one has the following (cf. \cite[Prop.~1.3.2]{vogel}).

\begin{proposition}\label{prop:vogel}
One may write
\begin{equation}\label{eq:shape r}
 r=\begin{cases}
    \prod_{i<j}[x_i,x_j]^{a_{ij}}\cdot r' & \text{if }p\neq2 \\ 
    \prod_{i=1}^n x_i^{2a_{ii}}\cdot\prod_{i<j}[x_i,x_j]^{a_{ij}}\cdot r' & \text{if }p=2 
   \end{cases} \qquad r'\in F_{(3)},
\end{equation}
with $0\leq a_{ij}<p$, and these numbers are uniquely determined by $r$.
Moreover, one has an isomorphism $\mathrm{tr}\colon H^2(G)\to\F_p$ such that $\mathrm{tr}(\chi_i\cup\chi_j)=-a_{ij}$,
with $a_{ij}$ as in \eqref{eq:shape r}.
\end{proposition}

In particular, from Proposition~\ref{prop:vogel} one deduces that if $\chi_h\cup\chi_k\neq0$ for some $1\leq h\leq k\leq d$, then 
\begin{equation}\label{eq:cup-product}
 \chi_i\cup\chi_j=\dfrac{a_{ij}}{a_{hk}}\chi_h\cup\chi_k, \qquad\text{for all }1\leq i\leq j\leq d.
\end{equation}

Let $G^{ab}=G/[G,G]$ be the abelianization of $G$.
Then $F^{ab}\simeq\Z_p^d$, and one may choose a basis $\mathcal{X}$ of $F$ such that $r\equiv x_1^q\bmod [F,F]$,
with $q$ a power of $p$ or $q=0$ if $r\in[F,F]$ --- i.e., one has an isomorphism of abelian pro-$p$ groups
\begin{equation}\label{eq:q}
 G^{ab}\simeq\Z_p^d\qquad\text{or}\qquad G^{ab}\simeq\Z_p/q\Z_p\times\Z_p^{d-1}.
\end{equation}
In particular, with such basis one has $a_{ii}=0$ for $i\geq2$, if $p=2$.
Henceforth, $\mathcal{X}$ will always denote a fixed basis of $F$ satisfying this condition.

The following two propositions are crucial for applying the theory exposed in Section~\ref{sec:3}.

\begin{proposition}\label{prop:onerel}
Let $G$ be a finitely generated one-relator pro-$p$ group, with minimal presentation \eqref{eq:presentation G}, basis 
$\mathcal{X}$, and defining relation $r$, such that $q\neq 2$. 
Then the following are equivalent:
\begin{itemize}
 \item[(i)] $H^\bullet(G)$ is quadratic;
 \item[(ii)] $r\notin F_{(3)}$, i.e., not all $a_{ij}$ are equal to 0.
\end{itemize}
Moreover, if the above conditions are satisfied then $G$ is mild.
\end{proposition}

\begin{proof} 
Since $q\neq2$, one has $a_{ii}=0$ for all $i\in\{1,\ldots,d\}$. 
Thus, by Proposition~\ref{prop:vogel} one has $\chi_i\cup\chi_i=0$ for all $i$, and moreover the initial form of $r$ is a Lie polynomial, namely, it is a combination of commutators $[X_i,X_j]$.

Suppose first that $H^\bullet(G)$ is quadratic, and let $\mathcal{B}=\{\chi_1,\ldots,\chi_d\}$ be a basis of $H^1(G)$ dual to
$\mathcal{X}$.
If $p$ is odd, then $H^2(G)$ is a quotient of $\Lambda_2(H^1(G))$ of dimension 1, so that there are $i,j$,
with $1\leq i<j\leq d$, such that $\chi_i\cup\chi_j\neq0$.
If $p=2$, then $H^2(G)$ is a quotient of $S_2(H^1(G))$ of dimension 1, so that there are $i,j$, with $1\leq i<j\leq d$,,
such that $\chi_i\cup\chi_j\neq0$. In both cases, Proposition~\ref{prop:vogel} implies that $a_{ij}\neq0$ 
as $\mathrm{tr}colon H^2(G)\to\F_p$ is an isomorphism, and this yields (ii).

Suppose now that $a_{ij}\neq0$ for some $1\leq i\leq j\leq d$.
Then \cite[Thm.~5.9--(i)]{jochen:massey} implies that $G$ is mild.
In particular, $H^n(G)=0$ for $n\geq3$ by Proposition~\ref{prop:mild}.
We claim that $H^\bullet(G)$ is quadratic.
Since $H^2(G)$ is the 1-dimensional vector space generated by $\chi_i\cup\chi_j$, the algebra $H^\bullet(G)$ is 1-generated.
Moreover, the fact that $H^3(G)$ is trivial is a consequence of the relations which hold in $H^2(G)$.
Indeed, let $1\leq h<k<l\leq d$ be any triplet of indices.
Then either $\chi_h\cup\chi_k=0$, or $\chi_h\cup\chi_k=b(\chi_i\cup\chi_j)$ for some $b\in\F_p^\times$.
In the former case one has $\chi_h\cup\chi_k\cup\chi_l=0$, whereas in the latter case one has 
\[
 \chi_h\cup\chi_k\cup\chi_l=b(\chi_i\cup\chi_j)\cup\chi_l.
\]
Then again either $\chi_j\cup\chi_l=0$, or $\chi_j\cup\chi_l=b'(\chi_i\cup\chi_j)$ for some $b'\in\F_p^\times$, so that 
\[
 \chi_h\cup\chi_k\cup\chi_l=bb'(\chi_i\cup\chi_i\cup\chi_j)=0,
\]
as $\chi_i\cup\chi_i=0$.
Therefore, the relations which hold in $H^2(G)$ imply that $H^3(G)=0$.
\end{proof}

\begin{proposition}\label{prop:onerel2}
Let $G$ be a finitely generated one-relator pro-2 group, with minimal presentation \eqref{eq:presentation G}, basis 
$\mathcal{X}$, and defining relation $r$, such that $q= 2$. 
\begin{itemize}
 \item[(i)] If $H^\bullet(G)$ is quadratic, then $r\notin F_{(3)}$, i.e., not all $a_{ij}$ are equal to 0.
 \item[(ii)] If $a_{ij}\neq0$ for some $1\leq i<j\leq d$, then $G$ is mild and $H^\bullet(G)$ is quadratic.
\end{itemize}
\end{proposition}

\begin{proof}
 The proof of statement~(i) is the same as the proof of the implication (i)$\Rightarrow$(ii) in Proposition~\ref{prop:onerel}.
 
 If $a_{ij}\neq0$ for some $1\leq i<j\leq d$, then the initial form of $r$ is the polynomial
 \[
  \rho=X_1^2+\sum_{1\leq h<k\leq d}a_{hk}[X_h,X_k]\in\F_2\langle X\rangle,\qquad a_{hk}\in\{0,1\},
 \]
with $a_{ij}=1$.
Therefore, one may choose an order on the set $X=\{X_1,\ldots,X_d\}$ such that the leading monomial of the homogeneous polynomial $\rho$ is $X_iX_j$. 
Since $i\neq j$, the sequence $\{X_iX_j\}$ is a {\sl combinatorially free} sequence of monomials of degree 2 (cf. \cite[Def.~3.1]{jochen:massey}), and therefore the sequence $\{\rho\}$ is strongly free by \cite[Thm.~3.5]{jochen:massey}, and $r$ yields a mild presentation of $G$.
In particular, Proposition~\ref{prop:mild} yields $H^n(G)=0$ for all $n\geq3$.

In order to prove that $H^\bullet(G)$ is quadratic if $a_{ij}\neq0$ for $1\leq i<j\leq d$, note that by Proposition~\ref{prop:vogel} one has $\chi_i\cup\chi_j\neq0$, and this cup-product generates $H^2(G)$.
Hence, if $\chi_h\cup\chi_k\neq0$ for some $1\leq h\leq k\leq d$, then Proposition~\ref{prop:vogel} implies
$\chi_h\cup\chi_k=\chi_i\cup\chi_j$ (in particular, $\chi_1\cup\chi_1=\chi_i\cup\chi_j$).
Therefore, the same argument as in the proof of implication (ii)$\Rightarrow$(i) of Proposition~\ref{prop:onerel} --- taking any triplet of indices $1\leq h\leq k\leq l\leq d$ (thus allowing equal indices) --- shows that $H^3(G)=0$ because of the relations which hold in $H^2(G)$, and this completes the proof of statement~(ii).
\end{proof}


\subsection{Demushkin pro-$p$ groups}\label{ssec:demushkin}

Here we describe briefly the example we keep in mind while dealing with one-relator pro-$p$ groups: Demushkin groups
(see \cite[\S~III.9]{nsw:cohm} and \cite[\S~5.2]{MQRTW} for further details).

A Demushkin group is a finitely generated one-relator pro-$p$ group $G$ such that the cup-product
induces a perfect pairing $H^1(G)\times H^1(G)\to\F_p$.
Equivalently, a finitely generated one-relator pro-$p$ group $G$ is Demushkin if and only if it has a presentation
\eqref{eq:presentation G} with defining relation $r$ such that one of the following holds:
 \begin{itemize}
 \item[(a)] $d$ is even and $r=x_1^{p^f}[x_1,x_2][x_3,x_4]\cdots[x_{d-1},x_d]$ for some $f\in\{1,2,\ldots\infty\}$
 such that $p^f\neq2$;
 \item[(b)] $d$ is even, $p=2$ and $r=x_1^{2+\alpha}[x_1,x_2]x_3^{2f}[x_3,x_4]\cdots[x_{d-1},x_d]$
 for some $f\in\{2,3,\ldots,\infty\}$ and $\alpha\in4\Z_4$;
 \item[(c)] $d$ is odd, $p=2$ and $r=x_1^2x_2^{2f}[x_2,x_3][x_4,x_5]\cdots[x_{d-1},x_d]$
 for some $f\in\{2,3,\ldots,\infty\}$
\end{itemize}
(cf. \cite{labute:demushkin} and \cite[Thm.~3.9.11 and Thm.~3,9,19]{nsw:cohm}).
In particular, the only finite Demushkin group is the cyclic group of order 2 (case (c) with d=1).
In this case, the $\F_2$-cohomology algebra $H^\bullet(G)$ is the ring of polynomials in one indeterminate with coefficients
in $\F_2$. 
Otherwise, $H^n(G)=0$ for $n\geq3$. 
In both cases, $H^\bullet(G)$ is quadratic.

Moreover, the graded group algebra $\gr_\bullet\FpG$ is isomorphic to the polynomial algebra 
$\FpX/(f)$, with $X=\{X_1,\ldots,X_d\}$, and $f$ a polynomial such that one of the following cases holds:
\begin{itemize}
 \item[(a)] $d$ is even and $f=[X_1,X_2]+[X_3,X_4]+\ldots+[X_{d-1},X_d]$;
 \item[(b)] $d$ is even, $p=2$ and $f=X_1^2+[X_1,X_2]+[X_3,X_4]+\ldots+[X_{d-1},X_d]$;
 \item[(c)] $d$ is odd, $p=2$ and $f=X_1^2+[X_2,X_3]+[X_4,X_5]+\ldots+[X_{d-1},X_d]$.
\end{itemize}
We call such a graded algebra $\FpX/(f)$ a {\sl Demushkin graded $\F_p$-algebra}.

\begin{remark} 
A quadratic algebra $A_\bullet$ has a single relation as above if, and only if, 
$\ext_{A_\bullet}^{2,2}(\F_p,\F_p)$ has dimension 1 
and the cup-product induces a non-degenerate alternating pairing 
$ \ext_{A_\bullet}^{1,1}(\F_p,\F_p)\times \ext_{A_\bullet}^{1,1}(\F_p,\F_p)\to \F_p$.
\end{remark}

Finally, one has the following
(cf. \cite[Thm.~5.2]{MQRTW}).

\begin{proposition}
 If $G$ is a Demushkin group, then $\gr_\bullet\FpG$ is isomorphic to the quadratic dual of $H^\bullet(G)$,
 and both algebras are Koszul.
\end{proposition}


\subsection{Cohomology}\label{ssec:H}

The isomorphism $\mathrm{tr}\colon H^2(G)\to\F_p$ induces a skewcommutative pairing
$\mathrm{tr}(\_\cup\_)\colon H^1(G)\times H^1(G)\to\F_p$.
If this pairing is perfect, then $G$ is a Demushkin group by definition.
Otherwise, let $V_2=H^1(G)^\perp$ be the radical of $H^1(G)$ with respect to the cup-product ---
i.e., 
\[V_2=H^1(G)^\perp=\{\chi\in H^1(G)\mid \chi\cup\psi=0\text{ for all }\psi\in H^1(G)\}.\]
Then $H^1(G)=V_1\oplus V_2$, so that the cup-product induces a perfect pairing
$V_1\times V_1\to \F_p$.
Set $q$ and $x_1$ as in \eqref{eq:q}.
Then \eqref{eq:shape r} yields
\begin{equation}\label{eq:shape r 3}
 r\equiv\begin{cases} \prod_{i<j}[x_i,x_j]^{a_{ij}}\mod F_{(3)}, & \text{if }q\neq2\\
    x_1^2\cdot\prod_{i<j}[x_i,x_j]^{a_{ij}}\mod F_{(3)}, & \text{if }q=2.   \end{cases}
\end{equation}

\begin{proposition}\label{prop:H}
Set $V_1$ and $V_2$ as above.
Then 
\begin{equation}\label{eq:case0 cohomology}
 H^\bullet(G)\simeq A_\bullet\sqcap Q(V_2,V_2^{\otimes2}),
\end{equation}
where $A_1=V_1$ and $A_2\simeq H^2(G)$, with the cup-product inducing a perfect pairing $A_1\times A_1\to\F_p$.
\end{proposition}

\begin{proof}
Let \eqref{eq:presentation G} be a minimal presentation of $G$.
If $q\neq2$, then the pairing induced by $\mathrm{tr}$ is alternating, so that $m=\dim(V_1)$ is even.
Hence, $V_1$ decomposes into a direct sum of hyperbolic planes (cf. \cite[Prop.~3.9.16]{nsw:cohm}).
Therefore, one may find a basis $\sB_1=\{\chi_1,\ldots,\chi_m\}$ of $V_1$
which completes to a basis $\sB$ of $H^1(G)$ such that 
\begin{equation}\label{eq:cupprod case1}
 \chi_1\cup\chi_2=\chi_3\cup\chi_4=\ldots=\chi_{n-1}\cup\chi_n=1
\end{equation}
and $\chi_i\cup\chi_j=0$ in any other case for $i\leq j$.

If $q=2$, let $\sB_1=\{\chi_1,\ldots,\chi_m\}$ be a basis of $V_1$ with $\chi_1$ dual to $x_1$.
Then by Proposition~\ref{prop:vogel} one has $\chi_1\cup\chi_1=1$, and $m$ can be both odd or even.
Thus, by \cite[Prop.~4]{labute:demushkin} we may choose the basis $\sB_1$
and complete it to a basis $\sB$ of $H^1(G)$ such that 
\begin{eqnarray*}
&&  \chi_1\cup\chi_2=\chi_3\cup\chi_4=\ldots=\chi_{m-1}\cup\chi_m=1,\qquad\text{if } 2\mid m,\\
&&  \chi_2\cup\chi_3=\chi_4\cup\chi_5=\ldots=\chi_{m-1}\cup\chi_m=1,\qquad\text{if } 2\nmid m.
\end{eqnarray*}
and $\chi_i\cup\chi_j=0$ in any other case for $i\leq j$.
\end{proof}


\begin{remark}\label{rem:boh} 
If $q\neq2$, then necessarily $\dim(V_1)\geq2$.
On the other hand, if $q=2$ then one may have $\dim(V_1)=1$. Then
$A_\bullet$ is isomorphic to the polynomial algebra in one indeterminate $\F_2[\chi_1]$.
This is the only case when $H^\bullet(G)$ is quadratic and $H^3(G)\neq0$.
\end{remark}

\subsection{The graded group algebra}\label{ssec:gr}

Let $G$ and $V_1,V_2\subseteq H^1(G)$ be as above.
First we deal with the case $q=2$ and $\dim(V_1)=1$, since by Proposition~\ref{prop:onerel} and Remark~\ref{rem:boh},
this is the only case with $G$ not mild.

\begin{proposition}\label{prop:caseq2}
Let $G$ be a finitely generated one-relator pro-2 group with $H^\bullet(G)$ quadratic and $\dim(V_1)=1$.
Then one has an isomorphism of graded $\F_2$-algebras
$\gr_\bullet\FpG=\F_2\langle X\rangle/( X_1^2)$, with $X=\{X_1,\ldots,X_d\}$, $d=\dd(G)$.
\end{proposition}

\begin{proof}
 By \eqref{eq:shape r 3} and Proposition~\ref{prop:vogel}, one has $r=x_1^2\cdot t$ with $t\in F_{(3)}$,
 with $\mathcal{X}=\{x_1,\ldots,x_d\}$ the basis of $F$.
Thus, the initial form of $r$ is $X_1^2\in F_2\langle X\rangle$, with $X=\{X_1,\ldots,X_d\}$.
We claim that $I(R)$ is the restricted ideal of $\FpX$ generated by $X_1^2$.

The subgroup $R\subseteq F$ is the (pro-2 closure of the) normal closure of the pro-$2$-cyclic group generated by $r$.
For $n\geq 2$, the non trivial elements of the subspace
\[
 \frac{(R\cap F_{(n)})F_{(n+1)}}{F_{(n+1)}}\leq\frac{F_{(n)}}{F_{(n+1)}}
\]
are the initial forms of all the elements of $R$ of degree $n$.
Such elements of $R$ are products of elements of the form $[y,r^{2^{m-1}}]$, with $m,s\geq0$ such that $n=2^m+s$, and $y\in F_{(s)}$; and also $y^{-1}r^{2^{m-1}}y$ in case $n$ is a power of 2, with $2^m=n$ and $y\in F$.
Commutator calculus and the properties \eqref{eq:propzass} yield
$$r^{2^{m-1}}=(x_1^2\cdot t)^{2^{m-1}}\equiv x_1^{2^m}\mod F_{(2^m+1)}$$
for all $m\geq1$, and consequently
\begin{equation}\label{eq:caseq2}
\begin{split}
\left[y,r^{2^{m-1}}\right] = \left[y,x^{2^{m}}t_m\right]\equiv \left[y,x_1^{2^m}\right] &\mod F_{(2^m+1+s)},\\
y^{-1}\cdot r^{2^{m-1}}\cdot y = \left(y^{-1}x_1^{2^m}y\right)\cdot \left(y^{-1}t_my\right) \equiv x_1^{2^m} &\mod F_{(2^m+1)},
\end{split}
\end{equation}
with $t_m\in F_{2^m+1}$, for $y\in F$ as above.
Therefore, the space $(R\cap F_{(n)})F_{(n+1)}/F_{(n+1)}$ --- viewed as subspace of the space of homogeneous polynomials of degree $n$ in $\F_2\langle X\rangle$ --- is generated by the polynomials $[\wp(X),X_1^{2^m}]$, with $\wp(X)$ running through all Lie polynomials in $F_2\langle X\rangle$ of degree $s$, with $m,s\geq0$ such that $n=2^m+s$; together with the monomial $X_1^n$ in case $n$ is a power of 2.

Hence, $I(R)$ is generated by the monomial $X_1^2$ as restricted ideal of the restricted $\F_2$-Lie algebra $L(X)$, and 
Proposition~\ref{prop:rest} yields $\gr_\bullet\FpG\simeq\gr_\bullet\F_p[F]/(X_1^2)$.
\end{proof}

On the other hand, if $\dim(V_1)\geq2$ one has the following.

\begin{proposition}\label{prop:grFG}
Set $V_1$ and $V_2$ as above, and assume that $\dim(V_1)\geq2$.
 Then the graded group $\F_p$-algebra of $G$ decomposes as free product
\begin{equation}\label{eq:grFG}
 \gr_\bullet\F_p[G]\simeq A_\bullet\sqcup T^\bullet(V_2^*) ,
\end{equation}
where $A_\bullet$ is a Demushkin quadratic $\F_p$-algebra with $A_1=V_1^*$.
\end{proposition}

\begin{proof}
Set $m=\dim(V_1)$, and let \eqref{eq:presentation G} be a minimal presentation of $G$.
Also, let $\mathcal{B}$ be a basis of $H^1(G)$ as in the proof of Proposition~\ref{prop:H}.

Let $\mathcal{S}=\{x_1,\ldots,x_m,y_1,\ldots,y_{d-m}\}\subseteq F$ be the basis dual to $\mathcal{B}$.
Then by \eqref{eq:cup-product} from \eqref{eq:shape r 3} one obtains
\begin{equation}\label{eq:rel case1}
 r\equiv [x_1,x_2]^a[x_3,x_4]^a\cdots[x_{m-1},x_m]^a\mod F_{(3)},
\end{equation}
for some $a\in\{1,\ldots,p-1\}$, if $q\neq2$, and 
\begin{equation}\label{eq:rel case2}
 r\equiv\left\{\begin{array}{c} x_1^2[x_1,x_2][x_3,x_4]\cdots[x_{m-1},x_m]\mod F_{(3)}\qquad\text{if }2\mid m, \\
 x_1^2[x_2,x_3][x_4,x_5]\cdots[x_{m-1},x_m]\mod F_{(3)}\qquad\text{if }2\nmid m, \end{array} \right. 
\end{equation}
if $q=2$.
Since $[x,y]^a\equiv[x^a,y]\equiv[x,y^a]\bmod F_{(3)}$, after a suitable change of basis $\mathcal{S}$ we may assume that $a=1$ in \eqref{eq:rel case1}.
Therefore, after identifying $\gr_\bullet\F_p[G]=\FpX$, whith $X=\{X_1,\ldots,X_d\}$,
the initial form of $r$ in $F_{(2)}/F_{(3)}$ is the homogeneous polynomial
\begin{equation}\label{eq:rel polynomial case1}
 \rho=[X_1,X_2]+[X_3,X_4]+\cdots+[X_{m-1},X_m]\in \F_p\langle X\rangle,
\end{equation}
if $q\neq2$, and 
\begin{equation}\label{eq:relation polynomial case2}
\rho=\left\{\begin{array}{c}X_1^2+[X_1,X_2]+[X_3,X_4]+\cdots+[X_{m-1},X_m]\qquad\text{if }2\mid m, \\
X_1^2+[X_2,X_3]+[X_4,X_5]+\ldots+[X_{m-1},X_m]\qquad\text{if }2\nmid m, \end{array} \right.
\end{equation}
if $q=2$.
Since $G$ is mild, Proposition~\ref{prop:mild} yields the claim.
\end{proof}


\subsection{Demushkin groups as quotients}

In \cite{wurf}, T.~W\"urfel proved that if a field $\K$ contains all roots of 1 of order a power of $p$
and its absolute Galois group $G_{\K}$ is a finitely generated one-relator pro-$p$ group, then one has a short exact sequence of pro-$p$ groups
\begin{equation}\label{eq:ses wurfel}
 \xymatrix{ 1\ar[r] & N \ar[r] & G_{\K}\ar[r] & G_{\K}/N\ar[r] & 1 }
\end{equation}
where $N$ is free and $G_{\K}/N$ is a Demushkin group, and moreover the inflation map
$\mathrm{inf}_{U,N}^2\colon H^2(U/N,\Z/p^s)\to H^2(U,\Z/p^s)$ 
is an isomorphism for every open subgroup $U\subseteq G_{\K}$ containing $N$ and every $s\geq1$.

If one considers finitely generated pro-$p$ groups whose closed subgroups have quadratic $\F_p$-cohomology (see Remark~\ref{rem:quad} below), one obtains the following.

\begin{proposition}\label{thm:onerel}
 Let $G$ be a finitely generated one-relator pro-$p$ group such that every closed subgroup of $G$ has quadratic $\F_p$-cohomology.
Then there exists a free closed subgroup $S\leq G$ and a short exact sequence of pro-$p$ groups
\begin{equation}\label{eq:ses wurfel 2}
 \xymatrix{ 1\ar[r] & N \ar[r] & G\ar[r] & G/N\ar[r] & 1 }
\end{equation}
where $N$ is the normal closure of $S$ in $G$, $G/N$ is a Demushkin group, and $N$ is free if $H^3(G)=0$.
Moreover, one has the isomorphisms of quadratic algebras
\begin{equation}\label{eq:intro}
\begin{split}
 H^\bullet(G) &\simeq H^\bullet(S)\sqcap H^\bullet(G/N),\\
\gr_\bullet\FpG &\simeq\gr_\bullet\F_p[S]\sqcup\gr_\bullet\F_p[G/N].
\end{split}
\end{equation}
\end{proposition}

\begin{remark}\label{rem:quad}
 The assumption that {\sl every closed subgroup} of $G$ has quadratic $\F_p$-cohomology might seem quite restrictive and unnatural.
 Still, by the Rost-Voevodsky Theorem every closed subgroup of the (maximal pro-$p$ quotient of the) absolute Galois group of a field containing a root of 1 of order $p$ has quadratic $\F_p$-cohomology (cf. Remark~\ref{rem:RV} below).
 Thus, this assumption is natural in view of the application of Proposition~\ref{thm:onerel} to the Galois-theoretic case
 (cf. Theorem~\ref{thm:Galois}).
\end{remark}

\begin{proof}[Proof of Proposition~\ref{thm:onerel}]
 Let $S$ be the closed subgroup of $G$ such that the restriction morphism $\res_{G,S}^1\colon H^1(G)\to H^2(S)$
induces an isomorphism $H^1(G)^\perp\simeq H^1(S)$.
In particular, one has $\kernel(\res_{G,S}^1)=V_1$.
Therefore, the commutative diagram
\[\xymatrix@R=0.8truecm{ H^1(G)\times H^1(G) \ar[rr]^-{\cup}\ar@<5ex>@{->>}[d]_{\res_{G,S}^1}\ar@<-5ex>@{->>}[d]_{\res_{G,S}^1}
&& H^2(G)\ar[d]_{\res_{G,S}^2} \\  H^1(S)\times H^1(S) \ar[rr]^-{\cup} && H^2(S) }\]
implies that the lower horizontal arrow is trivial and thus $H^2(S)=0$,
as $H^\bullet(S)$ is quadratic.
Consequently, $S$ is a free pro-$p$ group (cf. \cite[Prop.~3.5.17]{nsw:cohm}).

Let $N\subseteq G$ be the normal closure of $S$ in $G$, and set $\bar G=G/N$.
Since $H^1(N)^{\bar G}\simeq H^1(S)$, the exact sequence 
\[
 \xymatrix{ 0\ar[r] & H^1(\bar G)\ar[r]^-{\inf_{G,N}^1} & H^1(G)\ar[r]^-{\mathrm{res}_{G,N}^1} & H^1(N)^{\bar G}\ar[r] 
 & H^2(\bar G)\ar[r]^-{\inf_{G,N}^2} & H^2(G) }
\]
induced by the quotient $G/N$
implies that $H^1(\bar G)\simeq V_1$ and that the inflation map $\inf_{G,N}^2\colon H^2(\bar G)\to H^2(G)$
is a monomorphism (cf. \cite[Prop.~1.6.7]{nsw:cohm}).
Thus, in the commutative diagram
\[\xymatrix@R=0.8truecm{ H^1(\bar G)\times H^1(\bar G) \ar[rr]^-{\cup}\ar@<5ex>@{>->}[d]_{\inf_{G,N}^1}\ar@<-5ex>@{>->}[d]_{\inf_{G,N}^1}
&& H^2(\bar G)\ar@{>->}[d]_{\inf_{G,N}^2} \\  H^1(G)\times H^1(G) \ar[rr]^-{\cup} && H^2(G) }\]
the upper line is a non-degenerate pairing --- in particular, $\bar G$ is a one relator pro-$p$ group too.
Therefore, $\bar G$ is a Demushkin group (cf. \S~\ref{ssec:demushkin}).
Moreover, $\mathrm{inf}_{G,N}^2$ is an isomorphism, so that if $H^3(G)=0$ then \cite[Prop.~1]{wurf:cd2} implies that $N$ is free --- recall that $H^3(G)=0$ if, and only if, $\dim(V_1)=1$ (cf. Remark~\ref{rem:boh}).

Finally, \eqref{eq:intro} follows form Proposition~\ref{prop:H} and Proposition~\ref{prop:grFG}, since 
$V_1\simeq H^1(G/N)$ and $V_2\simeq H^1(S)$.
\end{proof}


\section{Absolute Galois groups of fields}
\label{sec:main}

Hereinafter $\K$ will denote a field containing a root of 1 of order $p$.
Moreover, $G_{\K}(p)$ will denote the maximal pro-$p$ quotient of the absolute Galois group of $\K$ --- namely,
$G_{\K}(p)$ is the maximal pro-$p$ Galois group (i.e., the Galois group of the maximal pro-$p$ extension) of $\K$.


\subsection{Maximal pro-$p$ Galois groups}
\label{ssec:proofs}

Let $\K^\times$ denote the multiplicative group of $\K$.
By Kummer theory one has an isomorphism $\K^\times/(\K^\times)^p\simeq H^1(G_{\K})$.
Moreover, note that if $p=2$ then $q=2$ (where $q$ is defined for $G_{\K}(2)$ as in \eqref{eq:q}) only if
$\sqrt{-1}\notin\K$.

The Rost-Voevodsky theorem has the following fundamental consequence (see, e.g., \cite[p.~222]{ido:miln}).

\begin{proposition}\label{thm:BK}
 The $\F_p$-cohomology ring $H^\bullet(G_{\K})$ of the absolute Galois group of $\K$ is quadratic.
 In particular, the epimorphism $G_{\K}\to G_{\K}(p)$ induces an isomorphism of graded algebras
 $  H^\bullet(G_{\K}(p))\simeq H^\bullet(G_{\K})$.
\end{proposition}

\begin{remark}\label{rem:RV}
 If $S$ is a closed subgroup of $G_{\K}$, respectively of $G_{\K}(p)$, then $S$ is the absolute Galois group $G_{\mathbb{L}}$, repsectively the maximal pro-$p$ Galois group $G_{\mathbb{L}}(p)$, of a suitable extension $\mathbb{L}/\K$, and obviously $\mathbb{L}$ contains a root of 1 of oder $p$ as well.
 Then by Proposition~\ref{thm:BK}, the $\F_p$-cohomology algebra $H^\bullet(S)$ is again quadratic.
\end{remark}

By Proposition~\ref{thm:BK} and \eqref{eq:H1H2}, $G_{\K}(p)$ is one-relator if, and only if, $H^2(G_{\K})$ has dimension 1.
Recall that the {\sl cohomological $p$-dimension} of a profinite group $G$ is the non-negative integer $\mathrm{cd}_p(G)$
defined by 
\[
 \mathrm{cd}_p(G)=\max\{n\geq0\mid H^n(G,M)\neq0\text{ for all $p$-torsion $G$-modules $M$}\}
\]
(cf. \cite[Def.~3.3.1]{nsw:cohm}).
If $G$ is a pro-$p$ group, then $\mathrm{cd}_p(G)=\mathrm{cd}(G)$.
Let $N$ denote the kernel of the epimorphism $G_{\K}\to G_{\K}(p)$.
Since $H^1(G_{\K}(p))\simeq H^1(G_{\K})$, the group $N$ is $p$-perfect, i.e., $H^1(N,\F_p)=0$, and hence $\mathrm{cd}_p(N)=0$.
Moreover, if $\mathrm{cd}(G_{\K}(p))<\infty$, then \cite[Prop.~3.3.8]{nsw:cohm} implies that $\mathrm{cd}_p(G_{\K})=\mathrm{cd}(G_{\K}(p))$.
Furthermore, if $\dim (H^1(G_\K))<\infty$, then $\mathrm{cd}(G_{\K}(p))$ (and hence $\mathrm{cd}_p(G_{\K})$) is finite (cf. \cite[Prop.~4.1]{cq:bk}).

On the other hand, the group algebras $\gr_\bullet\F_p[G_{\K}]$ and $\gr_\bullet\F_p[G_{\K}(p)]$
are related as follows (cf. \cite[Rem~1.4]{MQRTW} and \cite{shalev}).

\begin{proposition}\label{prop:shalev}
If $H^1(G_{\K})$ if finite the epimorphism $G_{\K}\to G_{\K}(p)$ induces an isomorphism of graded algebras
 $ \gr_\bullet\F_p[G_{\K}]\simeq\gr_\bullet\F_p[G_{\K}(p)] $.
\end{proposition}

From the results of Section~\ref{sec:4}, we may prove Theorem~\ref{thm:main}.

\begin{theorem}\label{thm:Galois}
Suppose that $H^1(G_{\K})$ is finite and $\dim H^2(G_{\K})=1$.
Then one has isomorphisms of quadratic algebras
\begin{equation}\label{eq:thm}
  H^\bullet(G_{\K})\simeq A_\bullet\sqcap Q(V_2,V_2^{\otimes2})\qquad\text{and}
 \qquad \gr_\bullet\F_p[G_{\K}]\simeq B_\bullet\sqcup T_\bullet(V_2^*),
\end{equation}
where $V_2=H^1(G_{\K})^\perp$ (with respect to the pairing induced by the cup-product),
$A_1\simeq H^1(G_{\K})/V_2$ and $B_\bullet$ a Demushkin algebra (cf. \S~\ref{ssec:demushkin}).

Moreover, there exists a closed subgroup $\tilde S\leq G_{\K}$ with $\mathrm{cd}_p(\tilde S)=1$ such that
\[\begin{split}
   A_\bullet\simeq H^\bullet(G_{\K}/N_{\tilde S})\qquad&\text{and}\qquad B_\bullet\simeq\gr_\bullet\F_p[G_{\K}/N_{\tilde S}]\\
   Q(V_2,V_2^{\otimes2})\simeq H^\bullet(\tilde S)\qquad&\text{and}\qquad T_\bullet(V_2^*)\simeq\gr_{\bullet}\F_p[\tilde S].
  \end{split}\]
--- here $N_{\tilde S}$ denotes the closed normal subgroup generated by $\tilde S$, and $\mathrm{cd}_p(N_{\tilde\S})=1$ as well,
unless $p=2$ and $\sqrt{-1}\notin \K$.
\end{theorem}

\begin{proof}
By Theorem~\ref{thm:BK} and Proposition~\ref{prop:shalev}, it is enough to show \eqref{eq:thm}
for $H^\bullet(G_{\K}(p))$ and $\gr_\bullet\F_p[G_{\K}(p)]$.
By hypothesis, $G_{\K}(p)$ is a finitely generated one-relator pro-$p$ group,
thus \eqref{eq:thm} follows from Theorem~\ref{thm:onerel}.

Now let $\tilde S\subseteq G_{\K}$ be the lift of $S$, with $S$ and $N$ as in Proposition~\ref{thm:onerel} for $G=G_{\K}(p)$.
I.e., for $\pi\colon G_{\K}\to G_{\K}(p)$ the canonical projection, one has $\kernel(\pi|_{\tilde S})=\kernel(\pi)$.
By Remark~\ref{rem:RV}, both $S$ and $\tilde S$ have quadratic $\F_p$-cohomology.
In fact, $S$ is the maximal pro-$p$ quotient of $\tilde S$, and thus $\mathrm{cd}_p(\tilde S)=\mathrm{cd}(S)$ (cf. \cite[Prop.~3.3.8]{nsw:cohm}).

Moreover, one has an isomorphism $G_{\K}/N_{\tilde S}\simeq G_{\K}(p)/N$,
as $N_{\tilde S}$ is the lift of $N$ and thus $N_{\tilde S}/\kernel(\pi)\simeq N$.
Hence, $G_{\K}/N_{\tilde S}$ is a Demushkin group, and $\mathrm{cd}_p(N_{\tilde S})=\mathrm{cd}(N)$.
Finally, $\mathrm{cd}(N)=1$ unless $\dim(H^1(G_{\K}(p))/N)=1$, and this case occurs inly if $p=q=2$, namely, only if
$\sqrt{-1}\notin\K$.

Therefore, one may deduce all the claims of the statement from Proposition~\ref{thm:onerel}.
\end{proof}

Note that Theorem~\ref{thm:Galois} (in fact already Proposition~\ref{thm:onerel}) implies the statement of W\"urfel's result, but for the bijectivity of the maps $\mathrm{inf}_{U,N}^2$.
Corollary~\ref{cor:main} follows from Theorem~\ref{thm:Galois} together with Example~\ref{ex:tensor}
and Proposition~\ref{prop:operations koszul}, as Demushkin algebras are Koszul (cf. \S~\ref{ssec:demushkin}).

\begin{example} 
Let $G=G_\circ\ast_{\hat p}S$ be the free product (in the category of pro-$p$ groups) of a Demushkin
 group $G_\circ$ with a finitely generated free pro-$p$ group $S$.
 Then $G/N\simeq G_\circ$ (with $N$ the normal closure of $S$), and one has 
 $$H^\bullet(G)\simeq H^\bullet(G_\circ)\sqcap H^\bullet(S)\qquad\text{and}\qquad 
\gr_\bullet\F_p[G]\simeq \gr_\bullet\F_p[G_\circ]\sqcup \gr_\bullet\F_p[S].$$
Such group is a pro-$p$ groups of elementary type (see next subsection).
\end{example}

\begin{example}\label{ex:KZ} 
For $p$ odd let $G$ be the pro-$p$ group with minimal presentation
\[G=\left\langle x_1,x_2,x_3\mid [x_1,x_2]=x_3^q\right\rangle,\]
with $q>1$ a power of $p$.
Such pro-$p$ group satisfies all the conditions in W\"urfel's theorem, as stated in \cite[Thm.2]{kz},
and by Proposition~\ref{prop:onerel} its $\F_p$-cohomology algebra is quadratic.
In particular, one has $H^\bullet(G)\simeq H^\bullet(S)\sqcap H^\bullet(\bar G)$, with $S=\langle x_3\rangle$
and $\bar G=G/N\simeq\Z_p^2$, with $N$ the normal closure of $S$, and 
\[\gr_\bullet\FpG\simeq\F_p[X_1,X_2]\sqcup\F_p[X_3].\]
Yet, the group $G$ is not realizable as maximal pro-$p$ Galois group of any $\K$, by \cite[Thm.~4.2 and Thm.~8.1]{eq:kummer}.
\end{example}


\subsection{Koszul algebras of elementary type}

Let $\mu_{p^\infty}$ denote the group of roots of unity of order a power of $p$ contained in the maximal pro-$p$
extension of $\K$.
The maximal pro-$p$ Galois group $G_{\K}(p)$ acts on $\mu_{p^\infty}$, fixing the roots of order $p$ (as they lie in $\K$).
Since the subgroup of $\mathrm{Aut}(\mu_{p^\infty})$ which fixes the roots of order $p$ is isomorphic
to the (multiplicative) group $1+p\Z_p=\{1+p\lambda\mid\lambda\in\Z_p\}$, one has a homomorphism of pro-$p$ groups
\[
 \theta_{\K}\colon G_{\K}(p)\longrightarrow 1+p\Z_p,
\]
called the cyclotomic character.

Following \cite[\S~3]{efrat:small}, we call a {\sl cyclotomic pro-$p$ pair} a pair $(G,\theta)$ consisting of a 
finitely generated pro-$p$ group $G$ and a homomorphism $\theta\colon G\to1+p\Z_p$ (the homomorphism $\theta$ is
also called an {\sl orientation} of $G$, cf.~\cite{cq:bk,qw:cyclotomic}).
A cyclotomic pro-$p$ pair is realizable arithmetically if there exists a field $\K$ such that $G\simeq G_{\K}(p)$
and $\theta$ coincides with the cyclotomic character.
The class of cyclotomic pro-$p$ pairs {\sl of elementary type} is the class of cyclotomic pro-$p$ pairs containing
\begin{itemize}
 \item[(a)] any pair $(F,\theta)$, with $F$ a finitely generated free pro-$p$ group and $\theta\colon F\to1+p\Z_p$
any orientation (including the trivial group with trivial orientation $\theta$);
\item[(b)] any pair $(G,\theta)$, with $G$ a Demushkin group and $\theta\colon1+p\Z_p$ as defined in 
\cite[Thm.~4]{labute:demushkin} (including the cyclic group of order 2 with the non-trivial orientation
$\theta\colon\Z/2\to\{\pm1\}$); 
\end{itemize}
and such that any pair contained in this class which is not of the type (a) or (b) may be obtained iterating the following to operations:
\begin{itemize}
 \item[(c)] by taking the free product $(G_1,\theta_1)\ast(G_2,\theta_2)$ of two cyclotomic pro-$p$ pairs of elementary type,
given by the free pro-$p$ product $G_1\ast_{\hat p}G_2$;
\item[(d)] by taking the the semi-direct product $(\Z_p\rtimes G,\theta\circ\pi)$ of a cyclotomic pro-$p$ pair $(G,\theta)$,
defined by $gzg^{-1}=\theta(g)\cdot z$ for all $z\in\Z_p$ and with $\pi\colon \Z_p\rtimes G\to G$ the canonical 
projection.
\end{itemize}

\begin{remark}
 \begin{itemize}
  \item[(i)] The $\F_p$-cohomology algebra of a pro-$p$ group of elementary type is always quadratic.
  \item[(ii)] A one-relator pro-$p$ group $G$ is of elementary type if, and only if, $G$ is isomorphic to the free pro-$p$ product $G_\circ\ast_{\hat p} S$ of a Demushkin group$G_\circ$ with a (possibly trivial) free pro-$p$ group $S$.
  In particular, by \cite[Thm.~12--(f)]{kz} the pro-$p$ group $G$ as in Example~\ref{ex:KZ} is not of elementary type.
 \end{itemize}
\end{remark}

Not all Demushkin groups are known to occur as $G_{\K}(p)$ for some field $\K$.
On the other hand, if one takes only Demushkin groups which occur as $G_{\K}(p)$ for some field $\K$ in item (b) above, then all cyclotomic pro-$p$ pairs of elementary type obtained are realizable as maximal pro-$p$ Galois groups.
I.~Efrat's {\sl Elementary type Conjecture} states that if a cyclotomic pro-$p$ pair is realizable arithmetically,
then it is of elementary type (cf. \cite{ido:ETC}, see also \cite[Question~4.8]{efrat:ETC} and \cite[\S~10]{marshall}).

In analogy with the class of cyclotomic pro-$p$ pairs of elementary type, we define the classes of
Koszul graded algebras of $G$-elementary type and of $H$-elementary type.

\begin{definition}\label{defi:pbwetc} 
The class of {\sl Koszul graded $\F_p$-algebras} of {\sl $G$-elementary type $\mathcal{KET}_G$} is the smallest class of quadratic $\F_p$-algebras such that, for $V$ be any finite $\F_p$-vector space,
 \begin{itemize}
  \item[(a)] the free algebra $T^\bullet(V)$ is in $\mathcal{KET}_{G}$;
  \item[(b)] any Demushkin algebra $A_\bullet$ is in $\mathcal{KET}_{G}$ (including the trivial $\F_2$-algebra on one generator $Q(\F_2,\F_2^{\otimes2})$);
 \end{itemize}
and such that, for $V$ any finite $\F_p$-vector space,
\begin{itemize}
 \item[(c)] if $A_\bullet,B_\bullet$ are in $\mathcal{KET}_{G}$, then also the free product $A_\bullet\sqcup B_\bullet$
 is in $\mathcal{KET}_{G}$;
 \item[(d)] if $A_\bullet$ is in $\mathcal{KET}_{G}$, then also the symmetric tensor product
 $A_\bullet\otimes S_\bullet(V)$ is in $\mathcal{KET}_{G}$.
\end{itemize}
 Dually, the class of {\sl Koszul graded $\F_p$-algebras} of {\sl $H$-elementary type} $\mathcal{KET}_{H}$
 is the smallest class of quadratic $\F_p$-algebras such that, for $V$ be any finite $\F_p$-vector space,
 \begin{itemize}
  \item[(a')] the trivial algebra $A_\bullet=Q(V,V^{\otimes2})$ is in $\mathcal{KET}_{H}$ (including the case $V=0$);
  \item[(b')] the quadratic dual $A_\bullet^!$ of any Demushkin algebra $A_\bullet$ is in $\mathcal{KET}_{H}$ 
  (including the polynomial algebra on one indeterminate $\F_2[X]$);
 \end{itemize}
and such that, for $V$ be any finite $\F_p$-vector space,
\begin{itemize}
 \item[(c')] if $A_\bullet,B_\bullet$ are in $\mathcal{KET}_{H}$, then also the direct product $A_\bullet\sqcap B_\bullet$
 is in $\mathcal{KET}_{H}$;
 \item[(d')] if $A_\bullet$ is in $\mathcal{KET}_{H}$, then also the skew-symmetric tensor product
 $A_\bullet\wedge \Lambda_\bullet(V)$ is in $\mathcal{KET}_{H}$.
\end{itemize}
\end{definition}

By Example~\ref{ex:tensor}, Proposition~\ref{prop:operations koszul} and \cite[\S~5.2]{MQRTW},
all algebras of $G$-elementary and $H$-elementary type are in fact Koszul.
Combining the restults obtained in \cite[\S~5]{MQRTW} together Proposition~\ref{thm:onerel},
one deduces the following.

\begin{proposition}\label{prop:ETCPBW}
 Let $\K$ be a field containing a root of unity of order $p$, and assume that $H^1(G_{\K})$ is finite.
\begin{itemize}
 \item[(i)] If $(G_{\K},\theta)$ is of elementary type, then $\gr_\bullet\F_p[G_{\K}]$ and $H^\bullet(G_{\K})$
are Koszul algebras of $G$-, respectively $H$-, elementary type.
 \item[(ii)] If $G_{\K}(p)$ is one-relator, then $\gr_\bullet\F_p[G_{\K}]$ and $H^\bullet(G_{\K})$
are Koszul algebras of $G$-, respectively $H$-, elementary type.
\end{itemize}
\end{proposition}

Therefore, in analogy with Efrat's conjecture, we formulate the following refinement of
Conjecture~\ref{conj}.

\begin{conj}\label{ques:ETC}
 Let $\K$ be a field containing a root of unity of order $p$, such that $H^1(G_{\K})$ is finite.
 \begin{itemize}
  \item[(i)] The graded group algebra $\gr_\bullet\F_p[G_{\K}]$ is a Koszul algebra of $G$-elementary type. 
\item[(ii)] The $\F_p$-cohomology ring $H^\bullet(G_{\K},\F_p)$ is a Koszul algebra of $H$-elementary type. 
 \end{itemize}
\end{conj}

By Proposition~\ref{prop:ETCPBW} a positive solution of the Elementary Type Conjecture would imply a positive answer to Conjecture~\ref{ques:ETC}.

\bibliography{onerel1}

\providecommand{\bysame}{\leavevmode\hbox to3em{\hrulefill}\thinspace}
\providecommand{\MR}{\relax\ifhmode\unskip\space\fi MR }
\providecommand{\MRhref}[2]{%
  \href{http://www.ams.org/mathscinet-getitem?mr=#1}{#2}
}
\providecommand{\href}[2]{#2}
\begin{thebibliography}{10}

\bibitem{cem}
S.K. Chebolu, I.~Efrat, and J.~Min\'{a}\v{c}, \emph{Quotients of absolute
  {G}alois groups which determine the entire {G}alois cohomology}, Math. Ann.
  \textbf{352} (2012), no.~1, 205--221. \MR{2885583}

\bibitem{CMQ}
S.K. Chebolu, J.~Min\'{a}\v{c}, and C.~Quadrelli, \emph{Detecting fast
  solvability of equations via small powerful {G}alois groups}, Trans. Amer.
  Math. Soc. \textbf{367} (2015), no.~12, 8439--8464. \MR{3403061}

\bibitem{ddsms}
J.D. Dixon, M.P.F. du~Sautoy, A.~Mann, and D.~Segal, \emph{Analytic pro-{$p$}
  groups}, second ed., Cambridge Studies in Advanced Mathematics, vol.~61,
  Cambridge University Press, Cambridge, 1999. \MR{1720368 (2000m:20039)}

\bibitem{ido:ETC}
I.~Efrat, \emph{Orderings, valuations, and free products of {G}alois groups},
  Sem. Structure Alg{\'e}briques Ordonn{\'e}es, {U}niv. {P}aris {VII}
  \textbf{54} (1995).

\bibitem{efrat:ETC}
\bysame, \emph{Pro-{$p$} {G}alois groups of algebraic extensions of {$\bold
  Q$}}, J. Number Theory \textbf{64} (1997), no.~1, 84--99. \MR{1450486
  (98i:11096)}

\bibitem{efrat:small}
\bysame, \emph{Small maximal pro-{$p$} {G}alois groups}, Manuscripta Math.
  \textbf{95} (1998), no.~2, 237--249. \MR{1603329 (99e:12005)}

\bibitem{ido:miln}
\bysame, \emph{Valuations, orderings, and {M}ilnor {$K$}-theory}, Mathematical
  Surveys and Monographs, vol. 124, American Mathematical Society, Providence,
  RI, 2006. \MR{2215492 (2007g:12006)}

\bibitem{eq:kummer}
I.~Efrat and C.~Quadrelli, \emph{The {K}ummerian property and maximal pro-{$p$}
  {G}alois groups}, J. Algebra \textbf{525} (2019), 284--310. \MR{3911645}

\bibitem{jochen:massey}
J.~G{\"a}rtner, \emph{Higher {M}assey products in the cohomology of mild
  pro-{$p$}-groups}, J. Algebra \textbf{422} (2015), 788--820. \MR{3272101}

\bibitem{jacobson}
N.~Jacobson, \emph{Lie algebras}, Dover Publications, Inc., New York, 1979,
  Republication of the 1962 original. \MR{559927}

\bibitem{KLM}
D.~Karagueuzian, J.P. Labute, and J.~Min\'a\v{c}, \emph{The {B}loch-{K}ato
  conjecture and {G}alois theory}, Ann. Sci. Math. Qu\'ebec \textbf{35} (2011),
  no.~1, 63--73. \MR{2848031}

\bibitem{kz}
D.H. Kochloukova and P.A. Zalesskii, \emph{Free-by-{D}emushkin pro-{$p$}
  groups}, Math. Z. \textbf{249} (2005), no.~4, 731--739. \MR{2126211
  (2005j:20030)}

\bibitem{labute:demushkin}
J.P. Labute, \emph{Classification of {D}emushkin groups}, Canad. J. Math.
  \textbf{19} (1967), 106--132. \MR{0210788 (35 \#1674)}

\bibitem{lodval}
J.-L. Loday and B.~Vallette, \emph{Algebraic operads}, Grundlehren der
  Mathematischen Wissenschaften [Fundamental Principles of Mathematical
  Sciences], vol. 346, Springer, Heidelberg, 2012. \MR{2954392}

\bibitem{marshall}
M.~Marshall, \emph{The elementary type conjecture in quadratic form theory},
  Algebraic and arithmetic theory of quadratic forms, Contemp. Math., vol. 344,
  Amer. Math. Soc., Providence, RI, 2004, pp.~275--293. \MR{2060204}

\bibitem{MQRTW}
J.~Min\'a\v{c}, F.~W. Pasini, C.~Quadrelli, and N.~D. T\^an, \emph{Koszul
  algebras and quadratic duals in galois cohomology}, preprint, available at
  {\tt arXiv:1808.01695}, 2018.

\bibitem{nsw:cohm}
J.~Neukirch, A.~Schmidt, and K.~Wingberg, \emph{Cohomology of number fields},
  second ed., Grundlehren der Mathematischen Wissenschaften [Fundamental
  Principles of Mathematical Sciences], vol. 323, Springer-Verlag, Berlin,
  2008. \MR{2392026 (2008m:11223)}

\bibitem{PoliPosi}
A.~Polishchuk and L.~Positselski, \emph{Quadratic algebras}, University Lecture
  Series, vol.~37, American Mathematical Society, Providence, RI, 2005.
  \MR{2177131 (2006f:16043)}

\bibitem{positselsky:koszul1}
L.~Positselski, \emph{Koszul property and {B}ogomolov's conjecture}, Int. Math.
  Res. Not. (2005), no.~31, 1901--1936. \MR{2171198 (2006h:19002)}

\bibitem{positselsky:koszul2}
\bysame, \emph{Galois cohomology of a number field is {K}oszul}, J. Number
  Theory \textbf{145} (2014), 126--152. \MR{3253297}

\bibitem{posivis:koszul}
L.~Positselski and A.~Vishik, \emph{Koszul duality and {G}alois cohomology},
  Math. Res. Lett. \textbf{2} (1995), no.~6, 771--781. \MR{1362968 (97b:12008)}

\bibitem{priddy}
S.B. Priddy, \emph{Koszul resolutions}, Trans. Amer. Math. Soc. \textbf{152}
  (1970), 39--60. \MR{0265437 (42 \#346)}

\bibitem{cq:bk}
C.~Quadrelli, \emph{Bloch-{K}ato pro-{$p$} groups and locally powerful groups},
  Forum Math. \textbf{26} (2014), no.~3, 793--814. \MR{3200350}

\bibitem{cq:qc}
\bysame, \emph{Finite quotients of {G}alois pro-{$p$} groups and rigid fields},
  Ann. Math. Qu\'e. \textbf{39} (2015), no.~1, 113--120. \MR{3374753}

\bibitem{QSV}
C.~Quadrelli, I.~Snopce, and M.~Vannacci, \emph{On pro-{$p$} groups with
  quadratic cohomology}, preprint, available at {\tt arXiv:1906.04789}, 2019.

\bibitem{qw:cyclotomic}
C.~Quadrelli and Th.S. Weigel, \emph{Profinite groups with a cyclotomic
  {$p$}-orientation}, preprint, available at {\tt arXiv:1811.02250}, 2018.

\bibitem{rost}
M.~Rost, \emph{Norm varieties and algebraic cobordism}, Proceedings of the
  {I}nternational {C}ongress of {M}athematicians, {V}ol. {II} ({B}eijing,
  2002), Higher Ed. Press, Beijing, 2002, pp.~77--85. \MR{1957022}

\bibitem{shalev}
A.~Shalev, \emph{Polynomial identities in graded group rings, restricted {L}ie
  algebras and {$p$}-adic analytic groups}, Trans. Amer. Math. Soc.
  \textbf{337} (1993), no.~1, 451--462. \MR{1093218}

\bibitem{BK}
V.~Voevodsky, \emph{On motivic cohomology with {$\bold Z/l$}-coefficients},
  Ann. of Math. (2) \textbf{174} (2011), no.~1, 401--438. \MR{2811603}

\bibitem{vogel}
D.~Vogel, \emph{Masseyprodukte in der {G}aloiskohomologie von
  {Z}ahlk\"{o}rpern}, Ph.D. thesis, University of Heidelberg, 2004.

\bibitem{weibel}
Ch. Weibel, \emph{2007 {T}rieste lectures on the proof of the {B}loch-{K}ato
  conjecture}, Some recent developments in algebraic {$K$}-theory, ICTP Lect.
  Notes, vol.~23, Abdus Salam Int. Cent. Theoret. Phys., Trieste, 2008,
  pp.~277--305. \MR{2509183}

\bibitem{weibel:art}
\bysame, \emph{Vladimir {V}oevodsky}, Notices Amer. Math. Soc. \textbf{66}
  (2019), no.~4, 526--533. \MR{3889528}

\bibitem{thomas:FP}
Th.S. Weigel, \emph{Graded {L}ie algebras of type {FP}}, Israel J. Math.
  \textbf{205} (2015), no.~1, 185--209. \MR{3314587}

\bibitem{thomas:proc}
\bysame, \emph{Koszul {L}ie algebras}, Lie algebras and related topics
  (M.~Avitabile, J.~Feldvoss, and Th. Weigel, eds.), Contemporary Mathematics,
  vol. 652, American Mathematical Society, Providence, RI, 2015, pp.~254--255.

\bibitem{wurf}
T.~W{\"u}rfel, \emph{A remark on the structure of absolute {G}alois groups},
  Proc. Amer. Math. Soc. \textbf{95} (1985), no.~3, 353--356. \MR{806069
  (87b:12006)}

\bibitem{wurf:cd2}
T.~W\"{u}rfel, \emph{Extensions of pro-{$p$} groups of cohomological dimension
  two}, Math. Proc. Cambridge Philos. Soc. \textbf{99} (1986), no.~2, 209--211.
  \MR{817662}

\end{thebibliography}

\bibliographystyle{amsplain}

\end{document}